\documentclass[11pt]{amsart}
\usepackage{amsfonts}
\usepackage{amssymb}
\usepackage{graphicx, enumerate}
\usepackage{pstricks}
\usepackage{amsmath}
\usepackage{amsxtra, hyperref}
\usepackage{mathrsfs, amsthm, xifthen, verbatim, mathrsfs}

\newtheorem{thm}{Theorem}[section]
\newtheorem{lemma}[thm]{Lemma}
\newtheorem{prop}[thm]{Proposition}
\newtheorem{crl}[thm]{Corollary}

\newtheorem{thA}{Theorem}

\theoremstyle{definition}
\newtheorem{dfn}[thm]{Definition}
\newtheorem{rem}[thm]{Remark}

\newcommand{\RR}{\mathbb{R}}

\newcommand{\NN}{\mathbb{N}}
\newcommand{\spn}{\text{Span}}
\newcommand{\N}[1]{\mathcal{N}(#1)}
\newcommand{\M}[1]{\mathcal{M}(#1)}
\newcommand{\X}[1]{\mathcal{X}(#1)}
\newcommand{\map}[3]{#1:#2\longrightarrow #3}
\newcommand\ddfrac[2]{\frac{\displaystyle #1}{\displaystyle #2}}

\begin{document}

\title[{The realizability problem as a special case of TMP}]{The realizability problem as a special case of the infinite-dimensional truncated moment problem}
\author[R.E. Curto and M. Infusino]{Ra\'ul E. Curto and Maria Infusino}
\address{R.E. Curto, Department of Mathematics, University of Iowa, Iowa City, 52246, USA}
\email{raul-curto@uiowa.edu}
\address{M. Infusino, Dipartimento di Matematica e Informatica, Universit\'{a} degli Studi di Cagliari, Palazzo delle Scienze, Via Ospedale 72, 09124 Cagliari}
\email{maria.infusino@unica.it}

\subjclass[2020]{Primary 44A60, 47A57, 60G55, 28C05; Secondary 46J05, 28E99}
\keywords{
truncated moment problem, point process, measure, integral representation, linear functional.
}

\maketitle

\begin{abstract}
The realizability problem is a well-known problem in the analysis of complex systems, which can be modeled as an infinite--dimensional moment problem. More precisely, as a truncated $K-$moment problem where $K$ is the space of all possible configurations of the components of the considered system. The power of this reformulation has been already exploited in \cite{KuLeSp11}, where necessary and sufficient conditions of Haviland type have been obtained for several instances of the realizability problem.
In this article we exploit this same reformulation to apply to the realizability problem the recent advances obtained in \cite{CGIK2022} for the truncated moment problem for linear functionals on general unital commutative algebras. This provides alternative proofs and sometimes extensions of several results in \cite{KuLeSp11}, allowing to finally embed them in the unified framework for the infinite-dimensional truncated moment problem presented in \cite{CGIK2022}.
\end{abstract}

\section{Introduction}

In several applications dealing with the analysis of complex systems, the underlying random distribution is often unknown or computationally intractable. Thus, a common strategy is to extract relevant information about the system from some selected characteristics rather than from the distribution itself. For example, in statistical mechanics the most interesting quantities characterizing many-particle systems can be deduced from the first two \emph{correlation functions}, which can be efficiently approximated via different methods such as Percus-Yevick schemes and hyper-netted chains (see, e.g., \cite{HMD87}). The challenging question is then to establish whether or not such estimated functions actually are the correlation functions of some random distribution on all possible configurations of the particles of the system. This problem is longstanding in the statistical mechanics literature where it is well-known as \emph{realizability problem} and its first rigorous formulation its due to Percus (see, e.g., \cite{P64}). Since the correlation functions can be defined as factorial moments of probability measures on the space of all possible configurations of the system components, which is an infinite-dimensional space, the realizability problem can be reformulated as an infinite-dimensional version of the classical truncated moment problem. 

Classically, the moment problem has been indeed investigated in finite dimensional settings. More precisely, given $d\in\NN$ and $K\subseteq\RR^d$ closed, the \emph{$d$-dimensional $K-$moment problem} asks whether a given sequence of real numbers is actually the sequence of power moments of a certain probability measure supported in $K$, or equivalently, whether a given linear functional on the real polynomials in $d$ variables can be represented as an integral w.r.t. a probability measure supported in $K$. The moment problem is addressed to as \emph{truncated} when the starting sequence is finite, respectively as \emph{full} when the starting sequence is infinite and so all moments are prescribed. 

The infinite-dimensionality can enter the moment problem in several ways, e.g., allowing measures supported in infinite-dimensional vector spaces, considering a starting sequence of functions rather than scalars, taking linear functionals defined on not necessarily finitely generated algebras, etc. All such versions belong to the class of \emph{infinite-dimensional moment problems}, whose study brings up new beautiful mathematical challenges as well as fascinating applications which have been triggering the interest in this problem from the early days of the moment problem theory until nowadays (see, e.g., \cite{BS71}, \cite[Chapter 5, Section 2]{BK88}, \cite{BY75}, \cite{H75}, \cite[Section 12.5]{Schm90}, \cite{GKM14},  \cite{IKR14}, \cite{AJK15}, \cite{GKM16}, \cite{Schm18}, \cite{GIKM}, \cite{IK2020}, \cite{IKKM-ArX} on the theoretical side, and, e.g., \cite{K67, E78, T02, TS06, M12, ZT20, M23} on the more applied one).

To identify the realizability problem as an infinite-dimensional truncated moment problem, we use in the following the interpretation given in \cite{KLS07, KuLeSp11} in which each configuration of the system components corresponds to a sum of Dirac measures concentrated at the points where the components lie. Thus, if we denote by $X$ the space of all the positions of the system components, each compact region of $X$ must contain only a finite number of components of each configuration and so each configuration $\eta$ is the sum of Dirac measures concentrated at the points of a sequence in $X$ without accumulation points, i.e., $\eta$ is a Radon measure on the space $X$. Hence, the set $\N X$ of all possible configurations of the system components is a subset of the vector space $\M X$ of all signed Radon measures on $X$. Prescribing the correlation functions of a probability measure on $\N X$ up to some order $n$ means specifying its factorial moments up to order $n$, which is equivalent to specifying its power moments up to order $n$, and so the realizability problem can be interpreted as a truncated $K-$moment problem with $K\subseteq \N X$. 

The full moment problem for point processes (i.e., probability measures on $\N X$) has been already considered in the seventies, e.g., in \cite{Le75a}  an analogue of Riesz-Haviland's theorem was proved for this instance, while in \cite{Krick73} a characterization of point processes was given via diagonal restrictions of their moments. Also solutions via positive semi-definiteness conditions have been formulated under additional analyticity bounds or boundedness constraints on the moments, e.g., in \cite{BeKoKuLy99},  \cite{KK02}, \cite{KoKuOl02}, \cite{IK2020}. As for the truncated case, much less is known and the available results solve specific instances of the problem for $K\subseteq\N X$ often motivated by instances of the realizability problem, see, e.g., \cite{US74}, \cite{CaKuLeSp06}, \cite{Kor05}, \cite{KLS07}, \cite{KuLeSp11}, \cite{MolLach15}, \cite{CIK16}, \cite{CGIK2022}. In particular, the power of the reformulation of the realizability problem as an infinite-dimensional truncated moment problem has been exploited in \cite{KLS07, KuLeSp11} to get a class of necessary and sufficient conditions for the existence of a representing measure on certain $K\subseteq\N X$ given only the first two correlation functions. Their results are based on a generalization of the Riesz-Markov theorem to unbounded continuous functions which exploits Daniell's theory of integration. 

In this manuscript, we instead exploit this reformulation to apply to the realizability problem the recent advances obtained in \cite[Theorems 3.7 and 3.8]{CGIK2022} for the truncated moment problem for linear functionals on general unital commutative algebras, retrieving and in some cases extending the main results in \cite{KuLeSp11}. The case of compactly supported measures on $\N X$ has been already considered in \cite[Section 6]{CGIK2022}, where the relative result in \cite{KuLeSp11} has been derived from \cite[Theorem~3.7]{CGIK2022}. In the following we complete the work started in \cite[Section 6]{CGIK2022}, by retrieving also the results in \cite{KuLeSp11}  for the case $K=\N X$ from the generalized truncated Riesz-Haviland theorem in \cite[Theorem 3.8]{CGIK2022} and so completely embedding \cite{KuLeSp11} in the unified approach to the infinite-dimensional truncated moment problem settled in~\cite{CGIK2022}.

\section{Notation and Preliminaries}

Let us start by recalling the above mentioned interpretation of the realizability problem as infinite-dimensional moment problem.
In the following we consider a Hausdorff locally compact space $X$ whose topology has a countable basis, and hence, $X$ is $\sigma$--compact and Polish (complete, separable, and metrizable). A possible configuration of the system components can be then modeled as a \emph{point configuration} $\eta$ on $X$, namely, a Radon measure $\eta$ on $X$ taking as values either a non-negative integer or infinity, i.e., $\eta=\sum_{i\in I}\delta_{x_i}$ where $\delta_{x_i}$ denotes the Dirac measure supported at $x_i$, $(x_i)_{i\in I}$ is such that $x_i\in X$ with $I$ either $\mathbb{N}$ or a finite subset of $\mathbb{N}$ and if $I=\mathbb{N}$ then the sequence $(x_i)_{i\in I}$ has no accumulation points in $X$. The requirement that $\eta$ is a Radon measure ensures that any compact set of $X$ contains only finitely many system components. The space $\N X$ of all point configurations is known as \emph{point configuration space on $X$} (see, e.g., \cite{KK02})
and is a closed subset of the space $\M X$ of all signed Radon measures supported in $X$ endowed with the so-called vague topology $\tau$, i.e., the weakest topology on $\M X$ such the map $\M X\to\RR$, $\nu\mapsto \int_Xf d\mu$ is continuous for all $f\in\mathcal{C}_c(X)$, where $\mathcal{C}_c(X)$ denotes the space of all continuous real-valued functions compactly supported in $X$ (see \cite[Lemma~4.4]{Kall}).

A \emph{point process on $X$} is a Radon probability measure on $\M X$ which is supported in $\N X$ and intuitively is a random distribution of points in $X$ such that, with probability one, any compact set contains only finitely many of these points (see, e.g. \cite{DJ2003} for an overview).

Given a point process $\mu$ on $X$, its $n^{th}$ correlation function $\rho^{(n)}_{\mu}$ is defined as the expected value of the $n^{th}$ factorial power $\eta^{\odot n}$ of an element $\eta\in \N X$, that is, the symmetric $\rho^{(n)}_{\mu}\in\M {X^n}$ such that 
$$\int_{X^n} f d\rho^{(n)}_{\mu}=\int_{\N X}\left(\int_{X^n} f d\eta^{\odot n}\right)d \mu, \quad \forall f\in \mathcal{C}_c(X^n)?$$
The realizability problem exactly asks the converse question: given $\rho_1\in \M X, \rho_2\in \M{X^2}$ symmetric and $K\subseteq \N X$, does there exists a point process $\mu$ on $X$ supported in $K$ such that $\rho_1=\rho^{(1)}_{\mu}$ and $\rho_2=\rho^{(2)}_{\mu}$?
Since the $n^{th}$ correlation function of $\mu$ in one-to-one correspondence with $n^{th}$ moment function of $\mu$, that is the expected value of the $n^{th}$ tensor power $\eta^{\otimes n}$ of an element $\eta\in \N X$, the realizability problem is equivalent to the following $K-$moment problem where $K\subseteq \N X$.

\begin{dfn}\label{realiz-moment}
Given $m_1\in \M X, m_2\in \M{X^2}$ symmetric and $K\subseteq \N X$, does there exists a point process $\mu$ on $X$ supported in $K$ such that for each $n\in\{1,2\}$ we have that
$$\int_{X^n} f d m_n=\int_{\N X}\left(\int_{X^n} f d\eta^{\otimes n}\right)d \mu, \quad \forall f\in \mathcal{C}_c(X^n)?$$
\end{dfn}

To apply the general results in \cite{CGIK2022} to the problem in Definition \ref{realiz-moment}, we need to reformulate it in terms of linear functionals on a unital commutative algebra. Indeed, we would like to see it as a special case of the following general truncated moment problem.

\begin{dfn}\label{B-TKMP}
Let $A$ be a unital commutative $\RR$--algebra and assume that its character space $\X{A}$ (i.e., the space of all unital algebra homomorphisms from $A$ to $\RR$) is non-empty. For each $a\in A$, define $\hat{a}(\alpha):=\alpha(a)$ for all $\alpha\in\X{A}$ and endow $\X{A}$ with the weakest topology $\tau_{\X{A}}$ making $\hat{a}$ continuous for all $a\in A$. 

Given a closed subset $K$ of $\X{A}$,  a linear subspace $B$ of $A$, and a linear functional $\map{L}{B}{\RR}$, the \emph{$B$--truncated $K-$moment problem} asks whether there exists a non-negative Radon measure $\nu$
whose support is contained in $K$ such that 
\[
	L(b)=\int\hat{b}~d\nu\quad\textrm{ for all } b\in B.
\]
When this representation exists, $\nu$ is called a $K$--representing measure for~$L$.
\end{dfn}

Clearly, if a $K-$representing measure exists then the linear functional $L$ is \emph{$K-$positive}, that is, $L(b)\geq 0$ for all $b\in \mathrm{\mathrm{Pos}}_B(K)$ where $$\mathrm{Pos}_{B}(K):=\{b\in B: \hat{b}\geq 0 \ \text{on}\ K\}.$$

The choice of $A$ needed to be able to take $K\subseteq \N X$ in Definition \ref{B-TKMP} has been already identified in \cite[Section 6]{CGIK2022} and is the space $\mathscr{P}$ of all polynomials of the following form
$$a(\eta):=\sum_{j=0}^N f_j \eta^{\otimes j}, \ N\in\NN_0,\  f_0\in\RR,\  f_j\in\mathcal{C}_c(X^j), \eta\in\M X,$$
where for any $n\in\NN$ and $\nu\in \M X$ we define:
\begin{itemize}
\item the (symmetric) $n^{th}$ power $\nu^{\otimes n}$  of $\nu$ as 
$$
\nu^{\otimes n}(dx_1, \ldots, dx_n):=\nu(dx_1)\cdots\nu(dx_n)
$$ 
\item 
$
f_n\nu^{\otimes n}:=\int_{X^n}f_n(x_1, \ldots, x_n)\nu^{\otimes n}(dx_1, \ldots, dx_n)
$ for any $f_n\in\mathcal{C}_c(X^n)$ 
 and $f_0\nu^{\otimes 0}:=f_0$ for any $f_0\in\RR$
 \item for any $n,m\in\NN_0$, $f_n\in\mathcal{C}_c(X^n)$ and $g_m\in\mathcal{C}_c(X^m)$
\begin{equation}\label{mult-cross}
(f_n\nu^{\otimes n})(g_m\nu^{\otimes m})=(f_n\otimes g_m)\nu^{\otimes (n+m)}
\end{equation}
\end{itemize}

In \cite[Proposition 6.1]{CGIK2022} the following embedding is proved, which ensures that for $A=\mathscr{P}$ in Definition \ref{B-TKMP} we can take $K\subseteq\N X$ as $\N X\subseteq \M X\hookrightarrow\X{\mathscr{P}}$.

\begin{prop}\label{prop-chara}
The space $(\M X, \tau)$ is topologically embedded in the character space $(\X{\mathscr{P}}, \tau_{\X{\mathscr{P}}})$, i.e., the following map is a homeomorphism onto its image:
$$\begin{array}{llll}
\phi:&\M X&\to& \X{\mathscr{P}}\\
\ & \nu &\mapsto & \phi(\nu),
\end{array}
$$
where for any $a(\eta):=\sum\limits_{j=0}^N f_j \eta^{\otimes j}\in \mathscr{P}$ we define
$ \phi(\nu)(a):=\sum\limits_{j=0}^N \int_{X^j}f_j d\nu^{\otimes j}.
$
\end{prop}

It is now clear that the problem in Definition \ref{realiz-moment} is nothing but a $\mathscr{P}^{(2)}-$truncated $K-$moment problem with $K\subseteq\N X$, where for any $N\in\mathbb{N}$ we denote by $\mathscr{P}^{(N)}$ the linear subspace of $\mathscr{P}$ consisting of all polynomials of degree at most $N$, i.e.,{\small
$$
\mathscr{P}^{(N)}:=\left\{a(\eta)=f_0+\sum_{j=1}^N  f_j\eta^{\otimes j}\in\mathscr{P}: \ f_0\in\RR, f_j\in \mathcal{C}_c(X^j) \ \text{for} \ j=1,\ldots,N \right\}.
$$}

\begin{rem}\label{corresp-dot-cross}
In \cite{KuLeSp11} the authors consider the factorial $n^{\textrm{th}}$ power $\eta^{\odot n}$ instead of the $n^{\textrm{th}}$ power $\eta^{\otimes n}$ we considered above. Denote by $\widetilde{\mathscr{P}}$ the set of polynomials defined by replacing $\eta^{\otimes n}$ with  $\eta^{\odot n}$ in the above definition of $\mathscr{P}$. \ Then as sets $\mathscr{P}=\widetilde{\mathscr{P}}$ and there is a bijective correspondence between $K$--positive linear functionals on $\mathscr{P}$ and $K$--positive linear functionals on $\widetilde{\mathscr{P}}$ for any $K\subseteq \N X$. Hence, solving the realizability problem is equivalent to solve the $\mathscr{P}^{(N)}$-truncated $K-$moment problem.
\end{rem}

When $K$ is compact, the $\mathscr{P}^{(N)}-$truncated $K-$moment problem can be easily solved thanks to \cite[Theorem 3.7]{CGIK2022}  (see  \cite[Corollary 6.2]{CGIK2022}, here restated for convenience).

\begin{crl}\label{risultato1}
Let $K\subseteq \M X$ be compact, $N\in\mathbb{N}$ and $L:\mathscr{P}^{(N)}\to \RR$ be linear. \ There exists a $K$--representing measure for $L$ if and only if $L$ is $K-$positive.
\end{crl}

In light of Remark~\ref{corresp-dot-cross}, it becomes clear that: \cite[Theorem 3.4]{KuLeSp11} follows from applying Corollary~\ref{risultato1} for $N=2$ and $K=\mathcal{N}_D(X):=\{\eta=\sum_i\delta_{x_i}\in \N X: d(x_i, x_j)>D,\ \forall i\neq j\}$ where $D>0$ and $d$ is a metric for the topology on $X$ so that $(X, d)$ is a complete metric space. In addition, 
\cite[Proposition 3.9]{KuLeSp11} (resp. \cite[Corollary 3.10]{KuLeSp11}) can be obtained from Corollary~\ref{risultato1} for $N=2$ and $K$ any compact subset of $\N X$ (resp. for $N=2$ and $K=\mathcal{N}^Q(X):=\{\eta\in \N X: \eta(X)=Q\}$ or $K=\mathcal{N}^{\leq Q}(X):=\{\eta\in \N X: \eta(X)\leq Q\}$ for $Q\in\mathbb{N}$).

In \cite{KuLeSp11} the authors also prove solubility criteria for the case $K=\N X$, namely \cite[Theorem~3.14]{KuLeSp11} for $X$ compact and  \cite[Theorem~3.17]{KuLeSp11} $X$ non-compact, which are not covered by Corollary \ref{risultato1}. In the following we are going to show that both these results can be derived from \cite[Theorem 3.8]{CGIK2022}, which we restate here for the readers' convenience.

\begin{thA}\label{thm3.8}
Let $A$ be a unital commutative $\RR$--algebra. Suppose $K\subseteq\X{A}$ is closed and non-compact, $B\subsetneq A$ is a linear subspace, and there exists $p\in A\setminus B$ such that $\hat{p}\ge1$ on $K$, $B_p:=\spn(B\cup\{p\})$ contains~$1$, $B_p$ generates $A$ and the following holds: 
\begin{equation} \label{p-existence2}
\textrm{for all } b \in B, \; \sup_{\alpha \in K} \left|\ddfrac{\hat{b}(\alpha)}{\hat{p}(\alpha)}\right| < \infty.      
\end{equation}

\noindent Let $\map{L}{B}{\RR}$ a $K$--positive linear functional. If $L$ has a $K$--positive extension $\bar{L}$ to $B_p$, then there exists a $K$--representing measure for $L$.
\end{thA}

Let us adapt to our setting some terminology introduced in \cite{KuLeSp11}:
\begin{itemize}
\item for $n\in\mathbb{N}_0$, a point process $\mu$ on $X$ is said to have \emph{finite local $n^{th}$ moment} if for every compact $\Lambda\subseteq X$ we have $\int_{\N X} (1\!\!1_{\Lambda}^{\otimes n}\eta^{\otimes n}) \mu(d \eta)<\infty$ 
\item for $n\in\mathbb{N}_0$, a point process $\mu$ on $X$ is said to have \emph{finite $n^{th}$ moment} if $\int_{\N X} (1\!\!1_X^{\otimes n}\eta^{\otimes n}) \mu(d \eta)<\infty$.
\item when $X$ is compact, a \emph{restricted cubic polynomial} is any $q_{f_0, f_1, f_2, f_3}\in\mathscr{P}^{(3)}$ of the form 
\begin{equation}\label{restricted-pol}
q_{f_0, f_1, f_2, f_3}(\eta)=f_0+f_1\eta+f_2\eta^{\otimes 2}+f_31\!\!1_X^{\otimes 3} \eta^{\otimes 3}, \quad f_0, f_3\in\RR, f_1\in \mathcal{C}_c(X),  f_2\in \mathcal{C}_c(X^2).
\end{equation}
We denote by $\mathcal{R}$ the collection of all restricted cubic polynomials defined in \eqref{restricted-pol} above, and by $\langle\mathcal{R}\rangle$ the real unital algebra generated by $\mathcal{R}$.
\item when $X$ is non-compact and $0<\Gamma\in\mathcal{C}_0(X)$, a \emph{$\Gamma-$restricted cubic polynomial} is any polynomial of the form: 
\begin{equation}\label{chi-restricted-pol}
q^{\Gamma}_{f_0, f_1, f_2, f_3}(\eta)=f_0+f_1\eta+f_2\eta^{\otimes 2}+f_3\Gamma^{\otimes 3} \eta^{\otimes 3}, \quad f_0, f_3\in\RR, f_1\in \mathcal{C}_c(X),  f_2\in \mathcal{C}_c(X^2).
\end{equation}
We denote by $\mathcal{R}_{\Gamma}$ the collection of all $\Gamma-$restricted cubic polynomials defined in \eqref{chi-restricted-pol} above, and by $\langle\mathcal{R}_{\Gamma}\rangle$ the real unital algebra generated by $\mathcal{R}_{\Gamma}$.
\end{itemize}

\begin{lemma}\label{iso-characters}
The space $\left(\X{\langle\mathcal{R}_{\Gamma}\rangle}, \tau_{_{\X{\langle\mathcal{R}_{\Gamma}\rangle}}}\right)$ is topologically isomorphic to $\X{\mathscr{P}}\times \mathbb{R}$ endowed with the product topology.
\end{lemma}
\proof
Set $\mathfrak{g}_{_{\Gamma}}(\eta):= \Gamma^{\otimes 3}\eta^{\otimes 3}$ for all $\eta\in\mathcal{M}(X)$ and observe that from the definition of $\Gamma-$restricted polynomials we immediately derive
$\mathcal{R}_{\Gamma} =\{p+\lambda \mathfrak{g}_{_{\Gamma}}: p\in\mathscr{P}^{(2)}, \lambda\in\mathbb{R}\} $. Then
\begin{eqnarray*}
\langle\mathcal{R}_{\Gamma}\rangle&=&\{p+\lambda \mathfrak{g}_{_{\Gamma}}^{\otimes n}: p\in\mathscr{P}, \lambda\in\mathbb{R}, n\in\mathbb{N}\}\\
&=&
\left\{\sum_{j=0}^Np_j\mathfrak{g}_{_{\Gamma}}^{\otimes j}: p_j\in\mathscr{P}, N\in\mathbb{N} \right\}=
\mathscr{P}[ \mathfrak{g}_{_{\Gamma}}]
\end{eqnarray*}

\noindent We define the map $$\begin{array}{ccc}
\Phi: \chi(\langle\mathcal{R}_{\Gamma}\rangle)&\longrightarrow& (\chi(\mathscr{P}) \times \mathbb{R})\\
\alpha &\mapsto& (\alpha\restriction_{\mathscr{P}}, \,\alpha(\mathfrak{g}_{_{\Gamma}}))
\end{array}$$
which is injective and continuous, in view of the chosen topologies. Also, for any $(\beta, \lambda)\in (\chi(\mathscr{P}) \times \mathbb{R})$ let us define
$$\alpha_{(\beta, \lambda)}\left(\sum_{j=0}^Np_j\mathfrak{g}_{_{\Gamma}}^{\otimes j}\right):=\sum_{j=0}^N\beta(p_j)\lambda^j;$$
clearly $\alpha_{(\beta, \lambda)}\in\chi(\langle\mathcal{R}_{\Gamma}\rangle)$. Then $\Phi$ is bijective and $\Phi^{-1}$, given by $\Phi^{-1}(\beta, \lambda)=\alpha_{(\beta, \lambda)}$, is also continuous.
\endproof 

\section{Main results}
In the following we are going to show that both \cite[Theorem~3.14]{KuLeSp11} and \cite[Theorem~3.17]{KuLeSp11} can be derived from Theorem~\ref{thm3.8}. Note that both of those results are for linear functionals on $\mathscr{P}^{(2)}$ but the reader will observe that the proof holds equally well for linear functionals on $\mathscr{P}^{(N)}$ for any $N>2$. 
\begin{thm}\label{appl-KLS}
Let $X$ be compact and $L:\mathscr{P}^{(2)} \to\RR$ be linear and $\N{X}-positive$. Then the following are equivalent:
\begin{enumerate}[a)]
\item there exists a $\N{X}-$representing measure for $L$ with finite third moment;
\item there exists $R>0$ such that
\begin{equation}\label{cond-KLS}
 \forall\  q_{f_0, f_1, f_2, f_3}\in \mathrm{Pos}_{\mathcal{R}}(\N{X}), \ L(f_0+f_1\eta+f_2\eta^{\otimes 2})+f_3R\geq 0.
\end{equation}
\end{enumerate}
\end{thm}

\proof\

\noindent(b) $\Rightarrow$ (a): Let us denote by $r(\eta):=1+1\!\!1_X^{\otimes 3} \eta^{\otimes 3}$ for all $\eta\in\M X$. Then 
$r(\eta)\geq 1$ for all $\eta\in\N{X}$ and $\spn(\mathscr{P}^{(2)}\cup \{r\})=\mathcal{R}$. Using Proposition~\ref{prop-chara}, we have the following embeddings $\N{X}\subseteq \M X\subseteq\X{\mathscr{P}}\subseteq \X{\langle \mathcal{R}\rangle}$. Also, $r\in \langle \mathcal{R}\rangle\setminus \mathscr{P}^{(2)}$ and, for all $b(\eta):=f_0+f_1\eta+f_2\eta^{\otimes 2}\in \mathscr{P}^{(2)}$, we obtain that:
\begin{eqnarray*}
 &\!\!&\sup\limits_{\sigma\in\N{X}}\left|\frac{b(\sigma)}{p(\sigma)}\right|\leq\sup\limits_{\sigma\in\N{X}}\frac{|f_0|+\int_X|f_1|d\sigma+\int_{X^2}|f_2|d\sigma^{\otimes 2}}{1+\sigma^{\otimes 3}(X)}\\
&\leq& \sup\limits_{\sigma\in\N{X}}\frac{\max\left\{|f_0|, \max\limits_{x\in X}|f_1(x)|, \max\limits_{x,y\in X}|f_2(x,y)|\right\}(1+\sigma(X)+\sigma(X)^2)}{1+\sigma^{\otimes 3}(X)}\\
&\leq& \max\left\{|f_0|, \max\limits_{x\in X}|f_1(x)|, \max\limits_{x,y\in X}|f_2(x,y)|\right\}\cdot\sup_{\sigma\in\N{X}}\left(\frac{1+\sigma(X)+\sigma(X)^2}{1+\sigma(X)^3}\right)=:\lambda_b.
\end{eqnarray*}
Since $X$ is compact and $\sigma$ is Radon, we have $\max\limits_{x\in X}|f_1(x)|<\infty$, $\max\limits_{x,y\in X}|f_2(x,y)|<\infty$ and $\sigma^{\otimes n}(X)=\sigma(X)^n<\infty$ for all $n\in\mathbb{N}$. Hence, $\max\left\{|f_0|, \max\limits_{x\in X}|f_1(x)|, \max\limits_{x,y\in X}|f_2(x,y)|\right\}<\infty$ and  $\sup\limits_{\sigma\in\N{X}}\left(\frac{1+\sigma(X)+\sigma(X)^2}{1+\sigma(X)^3}\right)=\sup\limits_{t\in\mathbb{R}}\left(\frac{1+t+t^2}{1+t^3}\right)<\infty$, which together ensure that $\lambda_b<\infty$. 

Since \eqref{cond-KLS} holds, the functional $L_R:\mathcal{R}\to\RR$ defined by 
$$L_R(q_{f_0, f_1, f_2, f_3}):=L(f_0+f_1\eta+f_2\eta^{\otimes 2})+f_3R,$$ for any $q_{f_0, f_1, f_2, f_3}(\eta)=(f_0+f_1\eta+f_2\eta^{\otimes 2}+f_31\!\!1_X^{\otimes 3} \eta^{\otimes 3})\in \mathcal{R}$,
is a $\N{X}-$positive linear extension of $L$ to $\mathcal{R}$. Hence, we can apply Theorem~\ref{thm3.8} for $A=\langle \mathcal{R}\rangle$, $K=\N{X}$, $B=\mathscr{P}^{(2)}$, $p=r$ and $\bar{L}=L_R$ (hence $B_p=\spn(B\cup\{p\})=\mathcal{R}$). This ensures that there exists a $\N{X}-$representing measure $\nu$ for $L$. Moreover, in Theorem~\ref{thm3.8}, it is shown that $\nu=\tilde{r}\mu$ where $\mu$ is a Radon measure on the compactification $\widetilde{\N X}$ of $\N X$ and $\tilde{r}$ is the continuous extension of $\frac{1}{r}$ to $\widetilde{\N X}$ and $\nu(\widetilde{\N X}\setminus \N X)=0$. Then $\nu$ has finite third moment since
$$\int_{_{\N X}} (1\!\!1_X^{\otimes 3}\eta^{\otimes 3}) \nu(d \eta)=\int_{_{\widetilde{\N X}}} 1\!\!1_X^{\otimes 3}\eta^{\otimes 3}\tilde{r}(\eta)\, \mu(d \eta)\leq \mu\left(\widetilde{\N X}\right)<\infty.
$$
Hence, (a) holds.\\

\noindent (a) $\Rightarrow$ (b): Since by (b) there exists a $\N{X}-$representing measure $\mu$ for $L$ with finite third moment $R:=\int_{\N{X}}1\!\!1_X^{\otimes 3} \eta^{\otimes 3}\mu(d\eta)$, for all $q_{f_0, f_1, f_2, f_3}\in \mathrm{Pos}_{\mathcal{R}}(\N{X})$ we have that
\begin{eqnarray*}
0&\leq& \int_{\N{X}} q_{f_0, f_1, f_2, f_3}(\eta)\mu(d\eta)\\
&=&\int_{\N{X}} (f_0+f_1\eta+f_2\eta^{\otimes 2}+f_31\!\!1_X^{\otimes 3} \eta^{\otimes 3}) \mu(d\eta)\\
&=&\int_{\N{X}} (f_0+f_1\eta+f_2\eta^{\otimes 2}) \mu(d\eta)+ f_3\int_{\N{X}}1\!\!1_X^{\otimes 3} \eta^{\otimes 3}\mu(d\eta)\\
&=&L(f_0+f_1\eta+f_2\eta^{\otimes 2})+f_3 R,
\end{eqnarray*}
i.e., \eqref{cond-KLS} is fulfilled. Hence, (b) holds.
\endproof

Using the correspondence in Remark \ref{corresp-dot-cross}, we immediately see that \cite[Theorem~3.14]{KuLeSp11} is nothing but Theorem~\ref{appl-KLS}, while \cite[Theorem~3.17]{KuLeSp11} corresponds to the following.

\begin{thm}\label{appl-KLS2}
Let $X$ be non-compact and $L:\mathscr{P}^{(2)}\to\RR$ linear and $\N X-$positive. Then the following are equivalent.
\begin{enumerate}[a)]
\item there exists a $\N{X}-$representing measure for $L$ with finite third local moments 
\item there exist $R>0$ and $0<\Gamma\in\mathcal{C}_0(X)$ such that
\begin{equation}\label{cond-KLS-chi}
 \forall \ q^{\Gamma}_{f_0, f_1, f_2, f_3}\in \mathrm{Pos}_{\mathcal{R}_{\Gamma}}(\N X),  L(f_0+f_1\eta+f_2\eta^{\otimes 2})+f_3R\geq 0  
\end{equation}
\end{enumerate}
\end{thm}

\proof\

\noindent(b) $\Rightarrow$ (a): Denoting by $r(\eta):=1+\Gamma^{\otimes 3} \eta^{\otimes 3}$ for all $\eta\in\M X$, we have that $r(\eta)\geq 1$ for all $\eta\in\N{X}$ and $\spn(\mathscr{P}^{(2)}\cup \{p\})=\mathcal{R}_{\Gamma}$. Combining Proposition~\ref{prop-chara} and Lemma \ref{iso-characters}, we have that $\N{X}$ is closed in $\X{\langle\mathcal{R}_{\Gamma}\rangle}$. Also $r\in \langle \mathcal{R}_{\Gamma}\rangle\setminus \mathscr{P}^{(2)}$ and for all $b(\eta):=f_0+f_1\eta+f_2\eta^{\otimes 2}\in \mathscr{P}^{(2)}$ we obtain that:
\begin{eqnarray*}
 \sup\limits_{\sigma\in\N{X}}\left|\frac{b(\sigma)}{p(\sigma)}\right|&\leq& \sup\limits_{\sigma\in\N{X}}\frac{|f_0|+\int_X|f_1|d\sigma+\int_{X^2}|f_2|d\sigma^{\otimes 2}}{1+\int_{X^3}\Gamma^{\otimes 3}d\sigma^{\otimes 3}}\\
&\leq& \sup\limits_{\sigma\in\N{X}}\frac{|f_0|+\int_X\lambda_1\Gamma d\sigma+\int_{X^2}\lambda_2\Gamma^{\otimes 2}d\sigma^{\otimes 2}}{1+\int_{X^3}\Gamma^{\otimes 3}d\sigma^{\otimes 3}}\\
&=& \sup\limits_{\sigma\in\N{X}}\frac{|f_0|+\lambda_1\int_X\Gamma d\sigma+\lambda_2\left(\int_{X}\Gamma d\sigma\right)^2}{1+\left(\int_{X}\Gamma d\sigma\right)^3}\\
&=& \sup\limits_{\sigma\in\N{X}}\frac{|f_0|+\lambda_1(\Gamma \sigma)+\lambda_2 (\Gamma \sigma)^2}{1+(\Gamma \sigma)^3}=:\lambda_b, 
\end{eqnarray*}
where $\lambda_1:=\sup_{x\in X} \frac{|f_1(x)|}{\Gamma(x)}<\infty$, $\lambda_2:=\sup_{(x,y)\in X^2}\frac{|f_2(x, y)|}{\Gamma(x)\Gamma(y)}<\infty$. Note that for any $\sigma\in\N{X}$ we have $(\Gamma \sigma)^n=\left(\int_X\Gamma d\sigma\right)^n<\infty$ for all $n\in\mathbb{N}$ and hence $\lambda_b=\sup\limits_{\sigma\in\N{X}}\left(\frac{|f_0|+\lambda_1(\Gamma \sigma)+\lambda_2 (\Gamma \sigma)^2}{1+(\Gamma \sigma)^3} \right)\leq\sup\limits_{t\in\mathbb{R}} \left(\frac{|f_0|+\lambda_1t+\lambda_2t^2}{1+t^3}\right)<\infty$.

Then, by the assumption \eqref{cond-KLS-chi}, the functional $L_R:\mathcal{R}_{\Gamma}\to\RR$ defined by 
$$L_R(q^{\Gamma}_{f_0, f_1, f_2, f_3}):=L(f_0+f_1\eta+f_2\eta^{\otimes 2})+f_3R$$ for any $q^{\Gamma}_{f_0, f_1, f_2, f_3}(\eta)=(f_0+f_1\eta+f_2\eta^{\otimes 2}+f_3\Gamma^{\otimes 3} \eta^{\otimes 3})\in \mathcal{R}_{\Gamma}$ is a $\N{X}-$positive linear extension of $L$ to $\mathcal{R}_{\Gamma}$. Hence, we can apply Theorem~\ref{thm3.8} for $A=\langle \mathcal{R}_{\Gamma}\rangle$, $K=\N{X}$, $B=\mathscr{P}^{(2)}$ and $p(\eta)=r(\eta)$ (hence $B_p=\spn(B\cup\{p\})=\mathcal{R}_\Gamma$). This ensures that there exists a $\N{X}-$representing measure $\nu$ for $L$. Moreover, in Theorem~\ref{thm3.8}, it is shown that $\nu=\tilde{r}\mu$ where $\mu$ is a Radon measure on the compactification $\widetilde{\N X}$ of $\N X$ and $\tilde{r}$ is the continuous extension of $\frac{1}{r}$ to $\widetilde{\N X}$ and $\nu(\widetilde{\N X}\setminus \N X)=0$. Then, for any $\Lambda\subset X$ compact, we have
$$\int_{\N X} (1\!\!1_\Lambda^{\otimes 3}\eta^{\otimes 3}) \nu(d \eta)=\int_{_{\widetilde{\N X}}} (1\!\!1_\Lambda^{\otimes 3}\eta^{\otimes 3})\tilde{r}(\eta) \mu(d \eta)\leq\mu(\widetilde{\N X})<\infty.$$
\smallskip

\noindent (a) $\Rightarrow$ (b):\

\noindent Since there exists a $\N{X}-$representing measure $\nu$ for $L$ with finite local third moments, by \cite[Lemma 1.16]{KuLeSp11}, there exists $0<\Gamma\in\mathcal{C}_0(X)$ such that $R:=\int_{\N X} (\Gamma^{\otimes 3}\eta^{\otimes 3}) \nu(d \eta)<\infty$. Then for all $q^\Gamma_{f_0, f_1, f_2, f_3}\in \mathrm{Pos}_{\mathcal{R}_{\Gamma}}(\N X)$ we have that
\begin{eqnarray*}
0&\leq& \int_{\N{X}} q^\Gamma_{f_0, f_1, f_2, f_3}(\eta)\nu(d\eta)\\
&=&\int_{\N{X}} (f_0+f_1\eta+f_2\eta^{\otimes 2}+f_3\Gamma^{\otimes 3} \eta^{\otimes 3}) \nu(d\eta)\\
&=&\int_{\N{X}} (f_0+f_1\eta+f_2\eta^{\otimes 2}) \mu(d\eta)+ f_3\int_{\N{X}}\Gamma^{\otimes 3} \eta^{\otimes 3}\nu(d\eta)\\
&=&L(f_0+f_1\eta+f_2\eta^{\otimes 2})+f_3 R,
\end{eqnarray*}
Hence, \eqref{cond-KLS-chi} is fulfilled.
\endproof

\begin{section}*{Acknowledgements}
The authors wish to thank the support to their collaboration received within the Programme Visiting Professor/Scientist 2022 of University of Cagliari financed by LR 7/2007 of the Autonomous Region of Sardinia. M. Infusino is a member of the GNAMPA group of INdAM.
\end{section}

\bibliographystyle{amsplain}

\end{document}